\documentclass{amsart}
\usepackage{graphicx}
\usepackage{fancyvrb}
\usepackage{amssymb}
\usepackage[pdfauthor={Richard J. Mathar},pdfkeywords={Mint Problem, Coins, Weighings},pdfsubject={Mathematics, Combinatorics, ApSimon, Mint},colorlinks=true,citecolor=blue]{hyperref}

\usepackage{amsmath}

\newtheorem{rem}{Remark}
\newtheorem{exa}{Example}

\begin{document}

\title{ApSimon's Mint Problem with Three or More Weighings}

\author{Richard J. Mathar} 
\email{mathar@mpia.de}
\urladdr{http://www.mpia.de/~mathar}
\address{Max-Planck Institute of Astronomy, K\"onigstuhl 17, 69117 Heidelberg, Germany}


\date{\today}

\begin{abstract}
ApSimon considered the problem of deciding
by a process of two weighings
on which of a known number of mints
emit either coins of a known genuine weight or emit coins of a different secondary but unknown weight.
The combinatorial problem consists of finding two sets
of coin numbers to be loaded on the tray for each of the weighings, and then to minimize
the total count of coins to be drawn from all mints for these two weighings.

This work yields numerical results for the generalized problem which allows
three or more weighings to settle which of the mints
produce either sort of coins.
\end{abstract}

\maketitle 

\section{Definitions} 
\subsection{Statement of the Problem}
Consider a set of $M$ mints issuing a coin with a known nominal weight $G$. There are
two suppliers for the coin material, each supplier producing the material for 
a fixed subset of the mints. Unfortunately one of the suppliers uses faulty material, so
for some of the mints all of their coins weigh $G(1+\epsilon)$ characterized
by some unknown nonzero excess $\epsilon$.
An investigator is equipped with an absolute scale, an allowance to draw
any number $C_m$ of coins from the mints numbered by $m$, $1\le m\le M$, and ordered to find out
which of the mints emit which of the two types of coins by weighing
two times a subset of these coins.

\subsection{Algebra of the Search Space}
The investigator's art of solving this problem is in finding a vector of
coin numbers $C_{1,m}\ge 0$ of the
first weighing and another vector of coin numbers $C_{2,m}\ge 0$ of the
second weighing
such that for any outcome of two measured weights a unique correspondence
exists to one of the $2^M$ variants of nominal and faulty coins of the mints.

ApSimon states the problem \cite{ApSimon}: what is the minimum number of coins
\begin{equation}
C(W,M)=\sum_{r=m}^M C_m
\end{equation}
involved in the $W=2$ weighings that allows to put either a label $d_m=0$ on mint $m$
if it produces the correct coins or a label $d_m=1$ on mint $m$ if it produces the faulty coins?

To rephrase, consider the first weight measured,
$\sum_{m=1}^M G C_{1,m}(1+d_m\epsilon)$,
and the second weight measured,
$\sum_{m=1}^M G C_{2,m}(1+d_m\epsilon)$ \cite{GuyAMM101}. The investigator may subtract
the known masses of the nominal coins, $\sum_m GC_{1,m}$ and $\sum_m GC_{2,m}$,
to reduce the two measurements to their excess weights
$\sum_{m=1}^M G C_{1,m}d_m\epsilon$ and
$\sum_{m=1}^M G C_{2,m}d_m\epsilon$. The unknown excess $\epsilon$
can be eliminated by considering
the measured known ratio of these two reduced weights because $\epsilon$ drops out,
\begin{equation}
X_1\equiv \frac{\sum_{m=1}^M C_{2,m}d_m}{ \sum_{m=1}^M C_{1,m}d_m}.
\end{equation}
The problem is solved if two vectors $C_{1,m}$ and $C_{2,m}$
are found such that all these ratios differ for the $2^M$ different
binary vectors $d_m$.

\begin{rem}
This could also be rephrased as mixing the coin numbers
such that no two of these excess vectors defined by plotting the points
of the first weighing and second weighing in a two-dimensional
coordinate system are collinear  \cite{AlonFCS37,KhovanovaArxiv1406}.
\end{rem}

The coins of any mint may be re-used for the second weighing.
If the total number $C=\sum_{m=1}^MC_m$ of the coins is sought to be minimal,
it would be wasteful not to use the full set $C_m$ of a mint with at
least one of the two weighings. So the set of coin numbers to be
searched for an optimum is evidently reduced to $0\le C_{w,m} \le C_m$
for $w=1,2$. Another obvious constraint is that from each mint $m$
at least one coin is to be put at the scale for at least one of the
weighings---otherwise no information of that $d_m$ would enter the weights.
So the cases $C_{1,m}=C_{2,m}=0$ do not need to be considered.

\section{Known Solutions for Two Weighings}
Guy and Nowakowski
found upper bounds of $C(2,6)\le 38$ for $M=6$ mints
and $C(2,7)\le 74$ coins for $M=7$ mints \cite{GuyAMM101}.
Li improved these upper bounds for two weighings to $31$ coins for 6 mints and 63 coins for 7 mints \cite{LiJTNU}.
Applegate settled the best value to $28$ coins for 6 mints, 51 for 7 mints and 90 coins for 8 mints \cite[A007673]{EIS}.

\begin{exa}
For $M=6$ mints the full information on the $d_m$ is extracted
by loading $C_{1,m}=(0,1,2,1,8,10)$ coins on the tray
for the first weighing and
$C_{2,m}=(1,2,2,5,5,0)$ coins for the second. This needs $C(2,6)=\sum_m C_m 
= \sum_m \max(C_{1,m},C_{2,m})=1+2+2+5+8+10=28$ coins
from all six mints for both weighings.
These two vectors of $C_{w,r}$ are not unique, because one could as well
combine $C_{1,m}=(0,2,1,1,8,10)$ and $C_{2,m}=(1,2,2,5,5,0)$ with the same total of $C(2,6)=28$
coins. There are two further solutions by just permuting the first and second weighing,
and there are further solutions by permuting
the enumeration of the $M$ mints, but the two solutions shown above are the only two
fundamentally different choices for the minimum of 28.
Even these two solutions are degenerate because permutation of the subset of the $C_{.,m}$
within a subset of constant $C_m$ (here $C_2=C_3=2$) does not cover different states of the $d_m$.
\end{exa}

\section{More than Two Weighings}
\subsection{Excess Weight Ratios}
Naturally the total number of coins needed becomes smaller if the
investigator may use a larger number $W$ of weighings, each with its own
set $C_{w,m}$ of coins, $1\le w\le W$. There is no new methodology to the
analysis but to require that the sets of ($W-1$) potentially measured
excess ratios
\begin{equation}
X_w\equiv \frac{\sum_{m=1}^M C_{w+1,m}d_m}{ \sum_{m=1}^M C_{w,m}d_m},\quad 1\le w< W,
\label{eq.rats}
\end{equation}
are $2^M$ different vectors of rational numbers as a function
of the binary state vectors $d_m$ \cite{BshoutyCOLT2009,HendyAMM87,GuyAMM102}.
\begin{rem}
Other definitions of the ratios may serve the same purpose. One might
for example use a constant reference value for $w$ in all the denominators,
or invert all ratios.
\end{rem}
Still 0-vectors of the form $C_{1,m}=C_{2,m}=\ldots =C_{w,m}=0$ do not
need to be considered because such an input cannot reveal information
on $d_m$. As for the case of two weighings,
minimization of the total number of coins requires
\begin{itemize}
\item
to use the full number $C_m$ in at least one
weighing to avoid waste,
\begin{equation}
C_m = \max(C_{1,m},C_{2,m},\ldots, C_{w,m}),
\end{equation}
\item
to use at least one coin of each mint in at least one weighing,
\begin{equation}
C_m \ge 1,
\end{equation}
\item
and to search for the minimum sum of coins purchased from all the mints,
\begin{equation}
C(W,M) = \min_{\{C_{w,m}\}} C
= \min_{\{C_{w,m}\}} \sum_{m=1}^M C_m.
\end{equation}
\end{itemize}

For the purpose of testing whether the ratios (\ref{eq.rats}) differ
for different sets of $d_m$, two different types of numbers are assigned
to the $X_w$ if the denominator is zero: If the numerator is positive,
$X_w=\infty$ as usual; if the numerator is also zero, a different
quantity $X_w=0/0$ is placed. Two different symbols for this case obviously
helps to reduce the number of coins needed, because a larger variation
of the components in the vectors $X_w$ helps to cover the $d_m$-space.

To illustrate this managing of zeros, consider the solution
\begin{equation}
C_{w,m}=\left(\begin{array}{cccc}
1 & 1 & 0 & 0\\
0 & 1 & 1 & 0 \\
0 & 0 & 1 & 2
\end{array}\right)
\label{eq.ex}
\end{equation}
to the problem with $W=3$ weighings and $M=4$ mints. Table \ref{tab.ex}
shows the $2^M$ different states, their $W$ excess weights and
ratios. For this
distribution of coins in the three weighings
the states of $d_m=(1,0,0,0)$ and $d_m=(1,0,0,1)$
can be distinguished, because the vector of the ratios is $(0,0/0)$
in the former case and differs from the vector of
ratios $(0,\infty)$  in the latter case. If $0/0$ and $\infty$ were
considered the \emph{same} ratio,
(\ref{eq.ex}) would not be flagged as a solution to the problem.

\begin{table}
\begin{tabular}{rrrr|ccc|cc}
$d_1$ & $d_2$ & $d_3$ & $d_4$ 
& $\sum_m C_{1,m}d_m$
& $\sum_m C_{2,m}d_m$
& $\sum_m C_{3,m}d_m$
& $X_1$
& $X_2$ \\
\hline
0 & 0 & 0 &  0 & 0 & 0 & 0 & $0/0$ & $0/0$ \\
0 & 0 & 0 &  1 & 0 & 0 & 2 & $0/0$ & $\infty$ \\
0 & 0 & 1 &  0 & 0 & 1 & 1 & $\infty$ & 1\\
0 & 0 & 1 &  1 & 0 & 1 & 3 & $\infty$ & 3\\
0 & 1 & 0 &  0 & 1 & 1 & 0 & 1 & 0\\
0 & 1 & 0 &  1 & 1 & 1 & 2 & 1 & 2\\
0 & 1 & 1 &  0 & 1 & 2 & 1 & 2 & 1/2\\
0 & 1 & 1 &  1 & 1 & 2 & 3 & 2 & 3/2\\
1 & 0 & 0 &  0 & 1 & 0 & 0 & 0 & 0/0\\
1 & 0 & 0 &  1 & 1 & 0 & 2 & 0 & $\infty$ \\
1 & 0 & 1 &  0 & 1 & 1 & 1 & 1 & 1\\
1 & 0 & 1 &  1 & 1 & 1 & 3 & 1 & 3\\
1 & 1 & 0 &  0 & 2 & 1 & 0 & 1/2 & 0\\
1 & 1 & 0 &  1 & 2 & 1 & 2 & 1/2 & 2\\
1 & 1 & 1 &  0 & 2 & 2 & 1 & 1 & 1/2\\
1 & 1 & 1 &  1 & 2 & 2 & 3 & 1 & 3/2\\
\end{tabular}
\caption{Decision table for 4 mints and 3 weighings with their reduced weights,
assuming coin counts specified by (\ref{eq.ex}).}
\label{tab.ex}
\end{table}

\subsection{Mints Equal to Weighings}
If the weighing number $W$ equals the number $M$ of mints, one could
use a single coin from a different mint $m$ for each of the weighings
and find individually one $d_m$ per weighing. Therefore a diagonal $C$-matrix
with column maximum $1$,
\begin{equation}
C_{w,m}=\delta_{w,m},\quad C_m=1
\end{equation}
suffices if $W\ge M$, and
\begin{equation}
C(W,M)=M,\quad W\ge M
\label{eq.diag}
\end{equation}
is an upper bound and also the optimum.

\subsection{Simple Bounds}
It is obvious that the number of coins needed is monotonous in both
variables:
\begin{itemize}
\item
\begin{equation}
C(W,M)\ge C(W+1,M)
\end{equation}
because
increasing the number of weighings does not require to increase
the number of coins to 
find the $d_m$. This is demonstrated by weighing
two times with the same assembly of coins, i.e., by duplicating a row
in the matrix $C_{w,m}$ of coins.
\item
\begin{equation}
C(W,M) \le C(W,M+1),
\end{equation}
because
increasing the number of mints requires no less coins to find
the $d_m$. This is proven by considering some minimizing solution $C_{w,m}$
with ratios $X_w$ for $M+1$ mints, chopping off
the component $d_{M+1}$ of the binary vector and removing the
associated ratios $X_w$  related to $d_{M+1}=1$, and observing
that in the reduced decision table all remaining $X_w$ vectors
are still pairwise different. (Example: delete all rows where $d_4=1$ and then the column $d_4$ in Table \ref{tab.ex},
which ends up in a decision table for 3 mints.)
\end{itemize}

\subsection{Solutions}
A list of one example of a matrix $C_{w,m}$ for the numerical solutions
that are found by exhaustive search follows. They have been computed with
a dedicated JAVA program reproduced in the \texttt{anc} directory.
\begin{multline}
C_{w,m}=\left(\begin{array}{cccc}
0&0&1&0\\
0&1&1&2\\
1&1&1&0\\
\end{array}\right)
\therefore C(3,4)=1+1+1+2=5.
\label{eq.C34}
\end{multline}
\begin{multline}
C_{w,m}=\left(\begin{array}{ccccc}
0&0&1&0&2\\
0&1&1&2&0\\
1&1&1&0&0\\
\end{array}\right)
\therefore C(3,5)=1+1+1+2+2=7.
\end{multline}
\begin{multline}
C_{w,m}=\left(\begin{array}{cccccc}
0&0&1&0&2&4\\
0&1&1&2&2&0\\
1&1&1&0&0&0\\
\end{array}\right)
\therefore C(3,6)=1+1+1+2+2+4=11.
\end{multline}
\begin{multline}
C_{w,m}=\left(\begin{array}{ccccccc}
0&0&1&3&0&2&4\\
1&1&1&3&1&1&0\\
0&1&1&0&3&3&4\\
\end{array}\right)
\\
\therefore C(3,7)=1+1+1+3+3+3+4=16.
\end{multline}
\begin{multline}
C_{w,m}=\left(\begin{array}{cccccccc}
1&0&1&0&0&5&0&5\\
0&0&1&2&4&1&2&5\\
0&1&1&2&2&5&5&0\\
\end{array}\right)
\\
\therefore C(3,8)=1+1+1+2+4+5+5+5=24.
\end{multline}
\begin{multline}
C_{w,m}=\left(\begin{array}{ccccc}
0&0&0&1&1\\
0&1&1&0&1\\
1&1&0&1&0\\
0&0&1&1&0\\
\end{array}\right)
\therefore C(4,5)=1+1+1+1+1=5.
\end{multline}
\begin{multline}
C_{w,m}=\left(\begin{array}{cccccc}
0&0&1&0&1&2\\
0&1&1&1&0&2\\
1&1&1&0&1&0\\
0&0&0&1&1&0\\
\end{array}\right)
\therefore C(4,6)=1+1+1+1+1+2=7.
\end{multline}
\begin{multline}
C_{w,m}=\left(\begin{array}{ccccccc}
0&0&1&0&1&2&1\\
0&1&1&1&0&2&3\\
1&1&1&0&1&0&0\\
0&0&0&1&1&0&3\\
\end{array}\right)
\\
\therefore C(4,7)=1+1+1+1+1+2+3=10.
\end{multline}
\begin{multline}
C_{w,m}=\left(\begin{array}{cccccccc}
0&0&1&0&0&2&2&2\\
0&1&1&1&2&0&2&2\\
1&1&1&1&0&2&2&0\\
0&0&0&1&2&2&2&0\\
\end{array}\right)
\\
\therefore C(4,8)=1+1+1+1+2+2+2+2=12.
\end{multline}
\begin{multline}
C_{w,m}=\left(\begin{array}{cccccc}
0&0&1&0&0&1\\
0&1&1&0&1&1\\
0&0&0&1&1&1\\
0&1&1&1&1&0\\
1&1&1&1&1&0\\
\end{array}\right)
\\
\therefore C(5,6)=1+1+1+1+1+1=6.
\end{multline}
\begin{multline}
C_{w,m}=\left(\begin{array}{ccccccc}
0&0&0&0&1&1&1\\
0&0&1&1&1&0&1\\
1&1&1&0&0&1&1\\
1&0&0&1&1&1&0\\
0&1&1&1&1&1&0\\
\end{array}\right)
\\
\therefore C(5,7)=1+1+1+1+1+1+1=7.
\end{multline}
\begin{multline}
C_{w,m}=\left(\begin{array}{cccccccc}
0&1&1&0&0&1&0&2\\
1&1&0&1&1&0&0&0\\
0&1&0&0&1&0&1&2\\
1&1&0&0&0&1&1&0\\
0&0&1&1&1&1&1&0\\
\end{array}\right)
\\
\therefore C(5,8)=1+1+1+1+1+1+1+2=9.
\end{multline}
\begin{multline}
C_{w,m}=\left(\begin{array}{ccccccccc}
0&1&1&0&0&1&0&2&0\\
1&1&0&1&1&0&0&0&0\\
0&1&0&0&1&0&1&2&3\\
1&1&0&0&0&1&1&0&3\\
0&0&1&1&1&1&1&0&1\\
\end{array}\right)
\\
\therefore C(5,9)=1+1+1+1+1+1+1+2+3=12.
\end{multline}
\begin{multline}
C_{w,m}=\left(\begin{array}{cccccccccc}
0&1&1&1&1&0&0&2&0&2\\
0&0&1&1&0&1&2&2&2&2\\
1&0&1&1&1&1&2&0&0&0\\
0&1&1&0&1&1&0&0&2&2\\
0&0&0&1&1&1&0&2&2&0\\
\end{array}\right)
\\
\therefore C(5,10)=1+1+1+1+1+1+2+2+2+2=14.
\label{eq.C510}
\end{multline}
\begin{multline}
C_{w,m}=\left(\begin{array}{ccccccccc}
0&0&0&1&0&1&0&1&1\\
0&1&1&0&1&1&0&0&1\\
0&1&0&1&1&0&1&0&1\\
1&1&1&1&1&0&0&1&1\\
1&1&0&0&0&1&1&1&1\\
0&0&1&1&1&1&1&1&0\\
\end{array}\right)
\\
\therefore C(6,9)=1+1+1+1+1+1+1+1+1=9.
\end{multline}
\begin{multline}
C_{w,m}=\left(\begin{array}{cccccccccc}
0&0&0&1&0&1&1&0&0&1\\
0&0&0&1&1&1&0&0&1&0\\
0&0&1&1&1&1&0&1&0&1\\
0&0&1&1&0&0&1&1&1&0\\
0&1&1&1&0&0&0&0&1&1\\
1&0&0&0&1&1&1&1&1&1\\
0&1&1&1&1&1&1&1&1&0\\
\end{array}\right)
\\
\therefore C(7,10)=1+1+1+1+1+1+1+1+1+1=10.
\label{eq.C710}
\end{multline}

\section{Summary}
Table \ref{tab.main}
shows the array $C(W,M)$  of the minimum number of  individual coins required with
$W$ weighings for $M$ mints
by collecting
result from equations (\ref{eq.C34})-(\ref{eq.C710}).
Entries below the diagonal are constant down the columns
according to (\ref{eq.diag}), and are not shown.
Entries with upper or lower bounds indicate that the
space of the $C_{w,m}$-matrices has not been scanned in full.
\begin{table}
\begin{tabular}{r|rrrrrrrrrr}
W$\backslash$ M &  1 & 2 & 3 & 4 & 5 &6 &7 &8 & 9& 10\\
\hline
2 & 1 & 2 & 4 & 8 & 15 & 28 & 51 & 90\\
3 & & & 3 & 5 & 7 & 11 & 16 & 24 & $\le 37$ & $\le 59$\\
4 & & &   & 4 & 5 & 7 & 10 & 12 & $>15, \le 18$ & $\le 28$ \\
5 & & & & & 5 & 6 & 7 &  9&  12 & 14\\
6 & & & & &   & 6 & 7 &  8& 9\\
7 & & & & &   &   & 7 &  8 & 9 & 10\\
\end{tabular}
\caption{The minimum number of coins $C(W,M)$ as a function of weighings $W$ and number of mints $M$.}
\label{tab.main}
\end{table}

\bibliographystyle{amsplain}
\bibliography{all}

\end{document}